\documentclass[12pt,epsf]{article}
\usepackage{amssymb,amsmath,eucal,amsthm}
\usepackage{graphicx}

\def\bce\{\begin{center}
\def\ece{\end{center}}

\newtheorem{proposition}{Proposition}[section]
\theoremstyle{definition}

\theoremstyle{remark}

\renewcommand{\proofname}{\bf Proof}

\title{Hopf bifurcation analysis  of pathogen-immune interaction
dynamics with delay kernel}

\author{ M. Neam\c tu$^{a}$\thanks{Corresponding author},   L. Buliga$^{b}$,   F.R. Horhat$^{c}$,    D. Opri\c s$^{b}$, A. T. Ceau\c su$^{b}$}
\date{ }

\begin{document}
\maketitle

\begin{tabular}{cccccccc}
\scriptsize{$^{a}$Department of Economic Informatics, Mathematics and Statistics, Faculty of Economics,}\\
\scriptsize{West University of Timi\c soara, str. Pestalozzi, nr. 16A, 300115, Timi\c soara, Romania,}\\
\scriptsize{E-mail:mihaela.neamtu@fse.uvt.ro,}\\
\scriptsize{$^{b}$ Department of Applied
Mathematics, Faculty of Mathematics,}\\
\scriptsize{West University of Timi\c soara, Bd. V. Parvan, nr. 4, 300223, Timi\c soara, Romania,}\\
\scriptsize{E-mails: opris@math.uvt.ro, larissa123us@yahoo.com}\\
\scriptsize{$^{b}$ Department of Biophysics and Medical Informatics,}\\
\scriptsize{University of Medicine and Pharmacy, Piata Eftimie Murgu, nr. 3, 300041, Timi\c soara, Romania,}\\
\scriptsize{E-mail: rhorhat@yahoo.com}\\

\end{tabular}
\begin{abstract}The aim of this paper is to study the steady states of the
ma\-the\-ma\-ti\-cal models with delay kernels which describe
pathogen-immune dynamics of many kinds of infectious diseases. In
the study of ma\-the\-ma\-ti\-cal models of infectious diseases it
is important to predict whether the infection disappears or the
pathogens persist. The delay kernel is described by the memory
function that reflects the influence of the past density of
pathogen in the blood and it is given by a nonnegative bounded
function $k$ defined on $[0,\infty )$ and is normated.
    By using the coefficient of  kernel $k$, as a bifurcation parameter, the
models are found to undergo a  sequence of Hopf bifurcation. The
direction and the stability criteria of bifurcation periodic
solutions are obtained by applying the normal form theory and the
center manifold theorems. Some numerical simulation examples for
justifying the theoretical results are also given.
\end{abstract}

\section{Introduction}

\hspace{0.5cm} The purpose of this paper is to study the Hopf
bifurcation of pathogen-immune dynamics in the steady states of
the mathematical model with delay kernel . Dynamical systems with
delay kernel have been studied for population dynamics and neural
networks [2].

We introduce a model which describes one of most known infectious
diseases, namely malaria infection. Our model is based on the
model from [6]. The model, as it was created, without delay, has
the feature that the interior equilibrium is always asymptotically
stable if it exists. For obtaining the natural behavior which was
observed experimentally, i.e. oscillatory be\-ha\-vior, it is
needed to make some adjustments to the model. One way is to
introduce some extra terms into the equations for a better
description of the interaction pathogens-immune system or other
way is to introduce the delay. Our contribution to the model lies
in introduction of the delay kernel, which is a natural thing to
do, because it is obvious for everyone that biological processes
do not take place instantaneous.

System (1) without the last term in the third equation, which
implies the effect of absorption of the pathogens into uninfected
cells, could be used for describing another well-known infectious
disease, the HIV infection, as was done in [3,7,10] or other
infectious diseases as hepatitis B virus infection [9] or
he\-pa\-ti\-tis C virus infection [7]. Because our model deals
with malaria infection in what follows we will say some words
about this disease. Malaria ranks high on the list of world health
problems by causing massive human and economic loss. The parasite
in this disease is called Plasmodium Falciparum and has a high
virulence which can cause even death. Malaria has some important
features that should be mentioned: first, the amount of variation
in disease severity observed in the field is remarkably high,
second, the immunity is virtually never sufficient to prevent
infection and third, transmission intensity in the field is highly
variable both temporally and geographically [4].

More precisely, consider the following system:

\begin{equation}
\begin{split}
&\dot{x}(t)=a_{1}-a_{2}x(t)-a_{3}x(t)\int\limits_{-\infty}^{0}
z(s)k(t-s)ds\\
 &\dot{y}(t)=-a_{4}y(t)+a_{3}x(t)\int\limits_{-\infty}^{0}z(s)k(t-s)ds\\
 &\dot{z}(t)=a_{4}a_{5}y(t)-a_{6}\int\limits_{-\infty}^{0}z(s)k(t-s)ds-a_{7}x(t)
\int\limits_{-\infty}^{0}z(s)k(t-s)ds,
\end{split}
\end{equation}
where $a_{i}$, $i=1,\cdots,7$ are positives constants and the
delay kernel, $k:[0,\infty)\rightarrow [0,\infty), \quad\mbox{k is
piecewise continuous},$ is assumed to  satisfy the following
properties:
\begin{eqnarray}\nonumber \\
\int\limits_{0}^{\infty}k(s)ds=1,\quad
\int\limits_{0}^{\infty}sk(s)ds< \infty.
\end{eqnarray}

It is also assumed that the system (1) is supplemented with
initial conditions of the form:
\begin{equation*}
\begin{split}
&x(0)=x^{\ast},\quad y(0)=y^{\ast},\quad
z(s)=\varphi_{1}(s),\quad s\in (-\infty,0], \\
&\varphi_{1}\quad\mbox{is bounded and continuous on $[0,\infty)$}.
\end{split}
\end{equation*}

The model contains three variables: the density of uninfected
cells x, the density of infected cells y and the density of
pathogens in blood z. Uninfected cells are recruited at a constant
rate $a_{1}$ from the source within the body, such as the bone
marrow and have the natural life expectancy  of $\frac{1}{a_{1}}$
days. Cells are infected by contact with pathogens and turn to
infected cells at rate $a_{3}\int\limits_{-\infty}^{0}
z(s)k(t-s)ds\nonumber$. Infected cells die at rate $a_{4}$. The
death of the cells results in the release of $a_{5}$ pathogens per
an infected cell and these pathogens have a life expectancy of
$\frac{1}{a_{6}}$. For $a_{7}=0$ we obtain the model of HIV
infection with delay kernel and for $a_{3}=a_{7}$ we obtain the
model of malaria infection with delay kernel.

In system (1), if the delay kernel  has the form
\begin{equation*}k(s)=\delta(s-\tau),\quad \tau\geq 0,
\end{equation*} where $\tau$ is a parameter which
denotes the effect of the past memories, then system (1) becomes:

\begin{equation}
\begin{split}
&\dot{x}(t)=a_{1}-a_{2}x(t)-a_{3}x(t)z(t-\tau) \\
&\dot{y}(t)=-a_{4}y(t)+a_{3}x(t)z(t-\tau)\\
&\dot{z}(t)=a_{4}a_{5}y(t)-a_{6}z(t-\tau)-a_{7}x(t)z(t-\tau)\\
&x(0)=x^{\ast},\quad y(0)=y^{\ast},\quad z(s)=\varphi_{1}(s),
\quad s\in [-\tau,0]
\end{split}
\end{equation}
For $\tau=0$ the system has been studied in [6].\\

 In system (1), if the kernel $k$ has the form
\begin{equation*}k(s)=qe^{-qs},
\end{equation*}called weak kernel, where q is a parameter varying
in $(0,\infty)$ which denotes the decay rate of the effect of the
post memories, then system (1) becomes:

\begin{equation}
\begin{split}
&\dot{x}(t)=a_{1}-a_{2}x(t)-a_{3}x(t)u(t)\\
&\dot{y}(t)=-a_{4}y(t)+a_{3}x(t)u(t)\\
&\dot{z}(t)=a_{4}a_{5}y(t)-a_{6}u(t)-a_{7}x(t)u(t)\\
&\dot{u}(t)=q(z(t)-u(t))
\\&x(0)=x^{\ast},\quad y(0)=y^{\ast},\quad z(s)=z^{\ast},\quad
u(0)=z^{\ast}
\end{split}
\end{equation}

In system (1), if the kernel $k$ has the form
\begin{equation*}k(s)=q^{2}se^{-qs},
\end{equation*}then system (1) becomes:

\begin{equation}
\begin{split}
&\dot{x}(t)=a_{1}-a_{2}x(t)-a_{3}x(t)v(t)\\
&\dot{y}(t)=-a_{4}y(t)+a_{3}x(t)v(t)\\
&\dot{z}(t)=a_{4}a_{5}y(t)-a_{6}v(t)-a_{7}x(t)v(t) \\
&\dot{u}(t)=q(z(t)-u(t))\\
&\dot{v}(t)=q(u(t)-v(t))\\
&x(0)=x^{\ast},\quad y(0)=y^{\ast},\quad z(0)=z^{\ast},\quad
u(0)=z^{\ast},\quad v(0)=z^{\ast}.
\end{split}
\end{equation}

This paper is organized as follows: in Section 2, the local
stability pro\-per\-ty and Hopf bifurcation of models (3), (4),
(5) are discussed and some sufficient conditions for stability are
derived. In Section 3, model (1) containing the general kernel is
further studied and both the direction and the stability of Hopf
bifurcation are analyzed by the normal form theory and the center
manifold theorem and some criteria for stability are derived.
Then, we consider two cases: in the first case $k$ is delta
function and in the second $k=qe^{-qs}$. Numerical simulations
will be shown, in order to justify the theoretical results.
Finally, conclusions are drawn with further research directions
given in Section 5.

\section{Local stability analysis and the Hopf bifurcation}

In this section, consider the local stability of the equilibrium
solution of system (1). From the special nature of the delay
kernel (2) embedded in system (1), we found out that an
equilibrium solution $(x_{0},y_{0},z_{0})$ of (1) is given by the
solution of the system:

\begin{equation}
\begin{split}
&a_{1}-a_{2}x-a_{3}xz=0\\
&a_{4}y-a_{3}xz=0\\
&a_{4}a_{5}y-a_{6}z-a_{7}xz=0.
\end{split}
\end{equation}

From (6) it results that if $0\leq a_{7}<a_{3}a_{5}$, $0<
a_{2}a_{6}<a_{1}(a_{3}a_{5}-a_{7})$, then system (6) has two
equilibria. The first one is $X_{1}=(\frac{a_{1}}{a_{2}},0,0)$ and
it represents the state in which the pathogens are absent. The
second is $X_{2}=(x_{0},y_{0},z_{0})$, where
\begin{equation*}
x_{0}=\frac{a_{6}}{a_{3}a_{5}-a_{7}},\quad
y_{0}=\frac{a_{1}(a_{3}a_{5}-a_{7})-a_{2}a_{6}}{a_{4}(a_{3}a_{5}-a_{7})},\quad
z_{0}=\frac{a_{1}(a_{3}a_{5}-a_{7})-a_{2}a_{6}}{a_{3}a_{6}}.
\end{equation*}
The equilibrium $X_{2}$ lies in the interior of the first
quadrant. Then we say that $X_{2}$ is an interior equilibrium and
represents the state in which the pathogens are present.\\
In what follows, the equilibrium $(x_{0},y_{0},z_{0})$ is
transformed to the origin, so the system (1) becomes:

\begin{equation}
\begin{split}
&\dot{u}_{1}(t)=-b_{1}u_{1}(t)-b_{2}\int\limits_{-\infty}^{0}k(-s)
u_{3}(t+s)ds-b_{3}u_{1}(t)\int\limits_{-\infty}^{0}k(-s)
u_{3}(t+s)ds\\
&\dot{u}_{2}(t)\!=\!b_{4}u_{1}(t)\!-\!b_{5}u_{2}(t)\!+\!b_{2}\int\limits_{-\infty}^{0}k(\!-\!s)
u_{3}(t\!+\!s)ds\!+\!b_3u_1(t)\!\!\int\limits_{-\infty}^{0}\!k(\!-\!s)
u_{3}(t\!+\!s)ds\\
&\dot{u}_{3}(t)\!=\!-\!b_{6}u_{1}(t)\!+\!b_{7}u_{2}(t)\!-\!b_{8}\!\!\int\limits_{-\infty}^{0}\!\!k(\!-\!s)
u_{3}(t\!+\!s)ds\!-\!b_{9}u_{1}(t)\!\!\int\limits_{-\infty}^{0}\!\!k(\!-\!s)
u_{3}(t\!+\!s)ds,
\end{split}
\end{equation}where

\begin{equation*}
\begin{split}
&b_{1}=a_{2}+a_{3}z_{0},\quad b_{2}=a_{3}x_{0}\quad
b_{3}=a_{3},\quad b_{4}=a_{3}z_{0},\quad b_{5}=a_{4}, \\
&b_{6}=a_{7}x_{0},\quad b_{7}=a_{4}a_{5},\quad
b_{8}=a_{6}+a_{7}x_{0},\quad b_{9}=a_{7}
\end{split}
\end{equation*}
 and
\begin{equation*}
u_{1}(t)=x(t)-x_{0},\quad
u_{2}(t)=y(t)-y_{0},\quad u_{3}(t)=z(t)-z_{0}.
\end{equation*}
Rewrite system (7) in the following matrix form:
\begin{equation*}\dot{u}(t)=Lu(t)+\int\limits_{-\infty}^{0}F(s)u(t+s)ds+H(u(t)),
\end{equation*}where
\begin{equation*}u(t)=(u_{1}(t),u_{2}(t),u_{3}(t))^{T},
\end{equation*}

\begin{equation}
\begin{split}
L=\left(\begin{array}{ccc}-b_{1}& 0 & 0\\
b_{4} & -b_{5} & 0\\-b_{6} & b_{7} & 0
\end{array} \right), F(s)=k(-s)\left(\begin{array}{ccc}0&0& -b_{2}\\0&0& b_{2}\\
0&0& -b_{8}\end{array}\right),\\ H(u(t))=\left(\begin{array}{ccc}
-b_{3}u_{1}(t)\int\limits_{-\infty}^{0}u_{3}(t+s)k(-s)ds\\
b_{3}u_{1}(t)\int\limits_{-\infty}^{0}u_{3}(t+s)k(-s)ds\\
-b_{9}u_{1}(t)\int\limits_{-\infty}^{0}u_{3}(t+s)k(-s)ds
\end{array}\right).
\end{split}
\end{equation}

The associated characteristic equation of the linearized system is
given by:
\begin{equation}\lambda^{3}+p_{2}\lambda^{2}+p_{1}\lambda+(r_{2}\lambda^{2}+r_{1}\lambda
+r_{0})\int\limits_{-\infty}^{0}k(-s)e^{\lambda s}ds=0,
\end{equation}where

\begin{equation*}
\begin{split}
&p_{2}=b_{1}+b_{5},\quad p_{1}=b_{1}b_{5},\quad r_{2}=b_{8},\quad
r_{1}=b_{8}(b_{1}+b_{5})-b_{2}b_{6}-b_{2}b_{7}\\
&r_{0}=b_{4}b_{7}-b_{1}b_{2}b_{7}-b_{1}b_{5}b_{8}-b_{2}b_{5}b_{6}.
\end{split}
\end{equation*}

\begin{proposition}If $k(s)=\delta(s-\tau)$, then
\begin{itemize}
\item[(i)]The characteristic equation (9) is given by
\begin{equation}\lambda^{3}+p_{2}\lambda^{2}+p_{1}\lambda+(r_{2}\lambda^{2}+r_{1}\lambda
+r_{0})e^{-\lambda\tau}=0;
\end{equation}
\item[(ii)]For $\tau=0$ the characteristic equation (10) is given
by
\begin{equation*}
\lambda^{3}+m_{2}\lambda^{2}+m_{1}\lambda+m_{0}=0,
\end{equation*}where
\begin{equation*}
m_{2}=p_{2}+r_{2},\quad m_{1}=p_{1}+r_{1},\quad m_{0}=r_{0};
\end{equation*}
\item[(iii)]If $\tau=0$, the equilibrium $X_{2}$ is asymptotically
stable if and only if
\begin{equation*}
m_{2}>0,\quad m_{1}>0,\quad m_{0}>0,\quad m_{1}m_{2}-m_{0}>0;
\end{equation*}
\item[(iv)]For $\tau=\tau_{0}$, given by
\begin{equation*}\tau_{0}=\frac{1}{\omega_{0}}\arctan\frac{r_{1}p_{2}\omega_{0}^{3}
-(\omega_{0}^{3}-\omega_{0}p_{1})(r_{0}-r_{2}\omega_{0}^{2})}{p_{2}\omega_{0}^{2}
(r_{0}-r_{2}\omega_{0}^{2})+r_{1}\omega_{0}(\omega_{0}^{3}-\omega_{0}p_{1})},
\end{equation*}
where $\omega_{0}$ is the positive root of the equation
\begin{equation}\nonumber
x^{6}+n_{1}x^{4}+n_{2}x^{2}+n_{3}=0,
\end{equation}
with
\begin{equation}\nonumber
n_{1}=p_{2}^{2}-2p_{1}-r_{2}^{2},\quad
n_{2}=p_{1}^{2}-r_{1}^{2}+2r_{0}r_{2},\quad n_{3}=-r_{0}^{2},
\end{equation}
there is a Hopf bifurcation.
\end{itemize}
\end{proposition}

\begin{proposition}If $k(s)=qe^{-qs}$, $s> 0$, $q>0$, then
\begin{itemize}
\item[(i)]The characteristic equation (9) is given by:
\begin{equation}\nonumber\lambda^{4}+(p_{2}+q)\lambda^{3}+(p_{1}+q(p_{2}+r_{2}))\lambda^{2}
+q(p_{1}+r_{1})\lambda+r_{0}q=0;
\end{equation}
\item[(ii)]The equilibrium $X_{2}$ is asymptotically stable if and
only if
\begin{equation}\nonumber
\begin{split}
D_{3}(q)&=((p_{1}+r_{1})(p_{2}+r_{2})-r_{0})q^{2}+
((p_{1}+r_{1})(p_{2}(p_{2}+\\&+r_{2})-r_{1})-2p_{2}r_{0})q+p_{2}(p_{1}-r_{0})>0;
\end{split}
\end{equation}
\item[(iii)]If there exists $q_{0}>0$ so that $D_{3}(q_{0})=0$ and
$\frac{dD_{3}(q)}{dq}|_{q=q_{0}}\neq 0$, then a Hopf bifurcation
occurs at $X_{2}$ as $q$ passes through $q_{0}$.
\end{itemize}
\end{proposition}

\section{Stability of the bifurcating periodic solutions: the general kernel case}

In this section, the stability of the bifurcating periodic
solutions of system (1) with the kernel satisfying (2) is studied.
For convenience, in the study of the Hopf bifurcation problem,
first we transform system (7) into an operator equation of the
form
\begin{equation}\nonumber\dot{u}_{t}=A(\mu)u_{t}+Ru_{t},
\end{equation}where
\begin{equation*}u=(u_{1},u_{2},u_{3})^{T},\quad
u_{t}=u(t+\theta),\quad \theta\in (-\infty,0),\quad \mu=a-a_{0}
\end{equation*}
and operators $A$ and $R$ are defined as

\begin{equation}\nonumber
A(\mu)\phi(\theta)=\left\{\begin{array}{rl}
\frac{d\phi(\theta)}{d\theta},&\theta\in (-\infty,0)\\
L\phi(\theta)+\int\limits_{-\infty}^{0}F(s)\phi(s)ds,&
\theta=0\end{array} \right.
\end{equation}

\begin{equation}\nonumber
R\phi(\theta)=\left\{\begin{array}{rl}(0,0,0)^{T},&\theta\in (-\infty,0)\\
(-b_{3}f_{1},b_{3}f_{1},-b_{9}f_{1})^{T},& \theta=0
\end{array}\right.
\end{equation}
where
\begin{equation}\nonumber
f_{1}=-\phi_{1}(0)\int\limits_{-\infty}^{0}k(-s)\phi_{3}(s)ds,
\end{equation} with $L, F$ defined in (8).\\
 As in [1], [5], the bifurcating periodic solutions $x(t,\mu)$ of (1)
 are indexed by a small parameter $\varepsilon\geq 0$. A solution
 $x(t,\mu(\varepsilon))$ has amplitude $O(\varepsilon)$, period
 $T(\varepsilon)$ and nonzero Floquet exponent $\beta(\varepsilon)$
 with $\beta(0)=0$. Under the present assumptions, $\mu, T$ and
 $\beta$ have expansions:

 \begin{equation}\nonumber
 \begin{split}&\mu=\mu_{2}\varepsilon^{2}+\mu_{4}\varepsilon^{4}+\cdots\\
 &T=\frac{2\pi}{\omega}(1+T_{2}\varepsilon^{2}+T_{4}\varepsilon^{4}+\cdots)\\
 &\beta=\beta_{2}\varepsilon^{2}+\beta_{4}\varepsilon^{4}+\cdots.
 \end{split}
 \end{equation}

 The sign of $\mu_{2}$ determines the direction of bifurcation,
 while $\beta_{2}$ determines the stability of
 $x(t,\mu(\varepsilon))$: asymptotically orbitally stable if
 $\beta_{2}<0$, but unstable if $\beta_{2}>0$.

 Next, the question of how to derive the coefficients in these
 expansions is addressed. For the applications from this paper, only $\mu_{2}, \tau_{2}$ and $\beta_{2}$ are
 computed here.

 We define
 the adjoint operator $A^{\ast}$ of $A$ as:
 \begin{equation}\nonumber
 A^{\ast}\psi(s)=\left\{\begin{array}{rl}
 -\frac{d\psi(s)}{ds},& s\in (0,\infty)\\
 L^{T}\psi(0)+\int\limits_{-\infty}^{0}F^{T}(s)\psi(-s)ds, & s=0,
 \end{array}\right.
 \end{equation}where $L^{T}$ and $F^{T}$ are transposes of
 matrices $L$ and $F$ respectively.

 Note that the operator $A$ depends on the bifurcation parameter
 $a$. According to Propositions 2.1, 2.2, Hopf bifurcation occurs when $a$
 passes through $a_{0}$. Let $\mu=a-a_{0}$. Then, Hopf
 bifurcation occurs when $\mu=0$. It is therefore reasonable to
 assume that $\varphi, \psi:[0,\infty)\rightarrow\mathbb{C}^{3}$.
 Define the bilinear form:

\begin{equation}\nonumber
<\phi,\psi>=\overline{\psi(0)}^{T}\phi(0)-\int\limits_{\theta=-\infty}
^{0}\int\limits_{\xi=0}^{\theta}\overline{\psi}^{T}(\xi-\theta)
F(\theta)\phi(\xi)d \xi d\theta.
\end{equation}
To determine the Poincare normal form of operator $A$, we need to
calculate the eigenvector $\phi$ of $A$ associated with eigenvalue
$\lambda_{1}=i\omega_{0}$ and the eigenvector $\phi^{\ast}$ of
$A^{\ast}$ associated with eigenvalue
$\lambda_{2}=\overline{\lambda_{1}}$.

\begin{proposition}
\begin{itemize}
\item[(i)]The eigenvector $\phi$ of $A$ associated with eigenvalue
$\lambda_{1}=i\omega_{0}$ is given by
$\phi(\theta)=ve^{\lambda_{1}\theta}$, $\theta\in(-\infty,0]$,
where $v=(v_{1},v_{2},v_{3})^{T}$ and
\begin{equation}v_{1}=-b_{2}(\lambda_{1}+b_{5})k^{1},
v_{2}=b_{2}(\lambda_{1}+b_{1}-b_{4})k^{1},
v_{3}=(\lambda_{1}+b_{1})(\lambda_{1}+b_{5}),
\end{equation}where
\begin{equation}\nonumber
k^{1}=\int\limits_{-\infty}^{0}k(-s)e^{\lambda_{1}s}ds;
\end{equation}
\item[(ii)] The eigenvector $\phi^{\ast}$ of $A^{\ast}$ associated
with eigenvalue $\lambda_{2}=\overline{\lambda_{1}}$ is given by
$\phi^{\ast}(s)=w e^{\lambda_{2}s}$, $s\in [0,\infty)$, where
$w=(w_{1},w_{2},w_{3})^{T}$ and

\begin{equation}
\begin{split}&w_{1}=\frac{b_{4}b_{7}-b_{6}(\lambda_{2}+b_{5})}{
b_{7}(\lambda_{2}+b_{1})\eta},\quad w_{2}=\frac{1}{\eta},\quad
w_{3}=\frac{\lambda_{2}+b_{5}}{b_{7}\eta}\\
 &\eta=\frac{b_{4}b_{7}-b_{6}(\lambda_{2}+b_{5})}{
b_{7}(\lambda_{2}+b_{1})}\overline{v}_{1}+\overline{v}_{2}+(\frac{\lambda_{2}+b_{5}}{b_{7}}-
(-b_{2}\frac{b_{4}b_{7}-b_{6}(\lambda_{2}+b_{5})}{b_{7}(\lambda_{2}+b_{1})}
\\&+b_{2}-b_{8}\frac{\lambda_{2}+b_{5}}{b_{7}})k^{(-1)})\overline{v}_{3},
\end{split}
\end{equation}
where
\begin{equation}\nonumber
k^{(-1)}=\int\limits_{-\infty}^{0}k(-s)e^{\lambda_{2}s}ds;
\end{equation}
\item[(iii)]We have:
\begin{equation}\nonumber<\phi^{\ast},\phi>=1,\quad
<\phi^{\ast},\overline{\phi}>=<\overline{\phi}^{\ast},\phi>=0,\quad
<\overline{\phi}^{\ast},\overline{\phi}>=1.
\end{equation}
\end{itemize}
\end{proposition}

Next, we construct the coordinates of the center of the manifold
$\Omega_{0}$ at $\mu=0$ ($a=a_{0}$) [1], [5]. Let

\begin{equation}\nonumber
\begin{split}&z(t)=<\phi^{\ast},x_{t}>\\
&w(t,\theta)=x_{t}- 2 Re \{z(t)\phi(\theta)\}.
\end{split}
\end{equation}

On the center  manifold $\Omega_{0}$,
$w(t,\theta)=w(z(t),\overline{z}(t),\theta)$, where
\begin{equation}\nonumber
w(z,\overline{z},\theta)=w_{20}(\theta)\frac{z^{2}}{2}+w_{11}(\theta)
z\overline{z}+w_{02}(\theta)\frac{\overline{z}^{2}}{2}+\cdots,
\end{equation}
$z$ and $\overline{z}$ are the local coordinates of the center
manifold $\Omega_{0}$ in the direction of $\phi$ and
$\phi^{\ast}$, respectively.\\
 For the solution $x_{t}\in\Omega_{0}$ of (1), notice that for
 $\mu=0$ we have:

\begin{equation}\nonumber\dot{z}(t)=\lambda_{1}z(t)+<\phi^{\ast},R(\omega+2Re\{z(t)\phi(\theta)\})>
\end{equation}
Rewrite this as
\begin{equation}\nonumber\dot{z}(t)=\lambda_{1}z(t)+g(z,\overline{z})
\end{equation}with
\begin{equation}\nonumber g(z,\overline{z})=\overline{\phi^{\ast}(0)}^{T}R(w(z,z,0)+2Re\{z(t)\phi(0)\})
\end{equation}

Further, expand the function $g(z,\overline{z})$ on the center
manifold $\Omega_{0}$ in powers of $z$ and $\overline{z}$:
\begin{equation}\nonumber
g(z,\overline{z})=g_{20}\frac{z^{2}}{2}+g_{11}z\overline{z}+g_{02}
\frac{\overline{z}^{2}}{2}+g_{21}\frac{z^{2}\overline{z}}{2}+\cdots
\end{equation}

\begin{proposition}For the system (1) we have:
\begin{itemize}
\item[(i)]
\begin{equation}
\begin{split}& g_{20}=-2b_{3}v_{1}v_{3}(\overline{w}_{1}-\overline{w_{2}}+
\frac{b_{9}}{b_{3}}\overline{w}_{3})k^{1}\\
 &g_{11}=-b_{3}(\overline{v}_{1}v_{3}+v_{1}\overline{v_{3}})
 (\overline{w}_{1}-\overline{w}_{2}+\frac{b_{9}}{b_{3}}\overline{w}_{3})
 (k^{1}+k^{(-1)})\\
 &g_{02}=-2b_{3}\overline{v}_{1}\overline{v_{3}}(\overline{w}_{1}-\overline{w}_{2}+
 \frac{b_{9}}{b_{3}}\overline{w}_{3})k^{(-1)};
\end{split}
\end{equation}
\item[(ii)]
\begin{equation}\nonumber
\begin{split}& w_{20}(\theta)=\frac{g_{20}}{\lambda_{1}}ve^{\lambda_{1}\theta}
-\frac{\overline{g_{20}}}{3\lambda_{1}}\overline{v}e^{\lambda_{2}\theta}+E_{1}e^{2\lambda_{1}\theta}\\
&w_{11}(\theta)=\frac{g_{11}}{\lambda_{1}}ve^{\lambda_{1}\theta}-
\frac{\overline{g_{11}}}{\lambda_{1}}\overline{v}e^{\lambda_{2}\theta}+E_{2},
\end{split}
\end{equation}where $E_{1}, E_{2}$ are the solutions of the
following system:
\begin{equation}\nonumber
\begin{split}&(A+k^{(2)}B-2\lambda_{1}I)E_{1}=b_{3}v_{1}v_{3}k^{1}(1,-1,-b_{9}/b_{3}^{2})^{T}\\
&(A+B)E_{2}=b_{3}(\overline{v}_{1}v_{3}+v_{1}\overline{v}_{3})(k^{1}+k^{(-1)})
(1,-1,-b_{9}/b_{3}^{2})^{T}\\
&k^{(1)}=\int\limits_{-\infty}^{0}k(-s)e^{\lambda_{1}s}ds,
k^{(2)}=\int\limits_{-\infty}^{0}k(-s)e^{2\lambda_{1}s}ds;
\end{split}
\end{equation}
\item[(iii)]
\begin{equation}
\begin{split}
&g_{21}=-2b_{3}(\overline{w}_{1}-\overline{w_{2}}+
\frac{b_{9}}{b_{3}}\overline{w}_{3})(v_{1}\int\limits_{-\infty}^{0}k(-s)
w_{311}(s)ds+\\&+\frac{1}{2}\overline{v}_{1}\int\limits_{-\infty}^{0}k(-s)
w_{320}(s)ds+v_{3}w_{111}(0)k^{1}+\frac{1}{2}\overline{v}_{3}w_{120}(0)k^{(-1)}),
\end{split}
\end{equation}
\end{itemize}
\end{proposition}
with $w_{20}(\theta)$=$(w_{120}(\theta)$, $w_{220}(\theta)$,
$w_{320}(\theta))^T$ and $w_{11}(\theta)$=$(w_{111}(\theta)$,
$w_{211}(\theta)$ $w_{311}(\theta))^T.$

Therefore, we will compute the following parameters:

\begin{equation}\nonumber
\begin{split}&c_{1}(0)=\frac{i}{2\omega_{0}}(g_{20}g_{11}-2|g_{11}|^{2}-\frac{1}{3}
|g_{02}|^{2})+\frac{g_{21}}{2}\\
&\mu_{2}=-\frac{Re\, c_{1}(0)}{Re \,\lambda'(a_{0})}\\
&T_{2}=-\frac{Im\,c_{1}(0)+\mu_{2}Im\,\lambda'(a_{0})
}{\omega_{0}}\\
&\beta_{2}=2Re\, c_{1}(0)\\
&T=\frac{2\pi}{\omega_{0}}(1+T_{2}\varepsilon^{2}+O(\varepsilon^{4})),\quad
\varepsilon^{2}=\frac{a-a_{0}}{\mu}+O(a-a_{0})^{2}.
\end{split}
\end{equation}

We have:

{\bf Theorem 3.1.} The sign of $\mu_{2}$ determines the directions
of the Hopf bifurcations: if $\mu_{2}>0(<0)$ the Hopf bifurcation
is supercritical (subcritical) and the bifurcating periodic
solutions exist for $a>a_{0}(<a_{0})$. The sign of $\beta_{2}$
determines the stability of the bifurcation periodic solutions.
They are both asymptotically orbitally  stable if $\beta_{2}< 0$,
but unstable if $\beta_{2}>0$. $T_{2}$ determines the period of
the bifurcating periodic
solutions: the period increases (decreases) if $T_{2}>0 (<0)$. \\

If $k(s)=\delta(s-\tau), \tau\geq 0$, then
$k^{1}=e^{\lambda_{1}\tau}$, $k^{(-1)}=e^{\lambda_{2}\tau}$,
$k^{2}=e^{2\lambda_{1}\tau}$ and $a_{0}=\tau_{0}$, where
$\tau_{0}$ is given by (20).
 In this case
 \begin{equation}\nonumber
 \begin{split}& \lambda'(\tau_{0})=
\frac{d\lambda}{d\tau}|_{\tau=\tau_{0},\lambda=\lambda_{1}}=\\&=
 \frac{\lambda_{1}(r_{2}\lambda_{1}^{2}+r_{1}\lambda_{1}+r_{0})}
 {(3\lambda_{1}^{2}+2p_{2}\lambda_{1}+p_{1})e^{\lambda_{1}\tau_{0}}+2r_{2}\lambda_{1}
 +r_{1}-(r_{2}\lambda_{1}^{2}+r_{1}\lambda_{1}+r_{0})\tau_{0}}
\end{split}
 \end{equation}

 If $k(s)=qe^{-qs}, s>0, q>0$, then
 $k^{1}=\frac{q}{\lambda_{1}+q}$,
 $k^{(-1)}=\frac{q}{\lambda_{2}+q}$and $a_{0}=q_{0}$, where
 $q_{0}$ satisfies $D_{3}(q_{0})=0$ ($D_{3}(q_{0})$ from Proposition 2.2), $
 \lambda_{1}=i\omega_{0}$,
 $\lambda_{2}=\overline{\lambda_{1}}$ and $\omega_{0}$ is given by
 \begin{equation}\nonumber
 \omega_{0}^2=\frac{q_{0}(p_1+r_{1})}{p_{2}+q_{0}}.
 \end{equation}
 In this case
 \begin{equation}\nonumber
\begin{split}&
 \lambda'(q_{0})=\frac{d\lambda}{dq}|_{q=q_{0},
 \lambda=\lambda_{1}}=\\&=-\frac{\lambda_{1}^{3}+(p_{2}+r_{2})\lambda_{1}^{2}+
 (p_{1}+r_{1})\lambda_{1}+r_{0}}{4\lambda_{1}^{3}+3(p_{2}+q_{0})\lambda_{1}^{2}
 +2(p_{1}+q_{0}(p_{2}+r_{2}))\lambda_{1}+q_{0}(p_{1}+r_{1})}.
\end{split}
 \end{equation}

From (11), (12), (13), (14) it results:

\begin{proposition} If $k(s)=qe^{-qs}$, $s>0$, $q>0$, then for system
(1) we have:
\begin{equation}\nonumber
\begin{split}
v_{1}=-\frac{b_{2}(\lambda_{1}+b_{5})q_{0}}{\lambda_{1}+q_{0}},\,
 v_{2}=\frac{b_{2}(\lambda_{1}+b_{1}-b_{4})q_{0}}{\lambda_{1}+q_{0}},\,
 v_{3}=(\lambda_{1}+b_{1})(\lambda_{1}+b_{5}),
\end{split}
 \end{equation}

\begin{equation}\nonumber
\begin{split}
 w_1=\frac{b_{4}b_{7}-b_{6}(\lambda_{2}+b_{5})}
 {b_{7}(\lambda_{2}+b_{1})\eta},\,
 w_2=\frac{1}{\eta},\,
 w_3=\frac{\lambda_{2}+b_{5}}{b_{7}\eta}
\end{split}
\end{equation}
\begin{equation}\nonumber
\begin{split}
 &\eta=\frac{b_{4}b_{7}-b_{6}(\lambda_{2}+b_{5})}
 {b_{7}(\lambda_{2}+b_{1})}\overline{v}_{1}+\overline{v}_{2}+(\frac{\lambda_{2}+b_{5}}{b_{7}}-
 (-b_{2}\frac{b_{4}b_{7}-b_{6}(\lambda_{2}+b_{5})}
 {b_{7}(\lambda_{2}+b_{1})}\\
 &+b_{2}-b_{8}\frac{\lambda_{2}+b_{5}}{b_{7}})
 \frac{q_{0}}{\lambda_{2}+q_{0}})\overline{v}_{3}\\
 \end{split}
\end{equation}
\begin{equation}\nonumber
\begin{split}
 &g_{20}=-2b_{3}v_{1}v_{3}(\overline{w}_{1}-\overline{w_{2}}+
\frac{b_{9}}{b_{3}}\overline{w}_{3}) \frac{q_{0}}{\lambda_{1}+q_{0}}\\
 &g_{11}=-b_{3}(\overline{v}_{1}v_{3}+v_{1}\overline{v_{3}})
 (\overline{w}_{1}-\overline{w}_{2}+\frac{b_{9}}{b_{3}}\overline{w}_{3})
 (\frac{q_{0}}{\lambda_{1}+q_{0}}+\frac{q_{0}}{\lambda_{2}+q_{0}})\\
 &g_{02}=-2b_{3}\overline{v}_{1}\overline{v_{3}}(\overline{w}_{1}-\overline{w}_{2}+
 \frac{b_{9}}{b_{3}}\overline{w}_{3})\frac{q_{0}}{\lambda_{2}+q_{0}}\\
 \end{split}
\end{equation}
\begin{equation}\nonumber
\begin{split}
 &g_{21}=-2b_{3}(\overline{w}_{1}-\overline{w}_{2}+\frac{b_{9}}
 {b_{3}}\overline{w}_{3})[v_{1}(\frac{g_{11}q_{0}}{\lambda_{1}
 (\lambda_{1}+q_{0})}v_{3}-\frac{\overline{g}_{11}q_{0}}{\lambda_{1}
 (\lambda_{2}+q_{0})}v_{3}+E_{32})\\&+\frac{1}{2}\overline{v}_{1}(
 \frac{g_{20}q_{0}}{\lambda_{1}(\lambda_{1}+q_{0})}v_{3}-
 \frac{\overline{g}_{20}q_{0}}{3\lambda_{1}
 (\lambda_{2}+q_{0})}\overline{v}_{3}+E_{31}\frac{q_{0}}{2\lambda_{1}+q_{0}})\\ &+
 v_{3}(\frac{g_{11}}{\lambda_{1}}v_{1}-\frac{\overline{g}_{11}}{\lambda_{1}}\overline{v}_{1}
 +E_{12})\frac{q_{0}}{\lambda_{1}+q_{0}}+\frac{1}{2}\overline{v}_{3}(\frac{g_{20}
 v_{1}}{\lambda_{1}}-\frac{\overline{g}_{20}\overline{v}_{1}}{3\lambda_{1}}+
 E_{11}\frac{q_{2}}{q_{2}+q_{0}})],
\end{split}
\end{equation}where $E_{1}=(E_{11},E_{21},E_{31})^{T}$, $E_{2}=(E_{12},E_{22},E_{32})^{T}$
satisfy the systems:
\begin{equation}\nonumber
\begin{split}&(A+\frac{q_{0}}{2\lambda_{1}+q_{0}}B-2\lambda_{1}I)E_{1}=
b_{3}v_{1}v_{3}\frac{q_{0}}{2\lambda_{1}+q_{0}}(1,-1,-\frac{b_{9}}{b_{3}^{2}})^{T} \\
&(A+B)E_{2}=b_{3}(\overline{v}_{1}v_{3}+v_{1}\overline{v}_{3})q_{0}
(\frac{1}{\lambda_{1}+q_{0}}+\frac{1}{\lambda_{2}+q_{0}})(1,-1,-\frac{b_{9}}{b_{3}^{2}})^{T}.
\end{split}
\end{equation}
\end{proposition}

\section{Numerical simulations}
For numerical simulations, we use Maple 9.5. We consider system
(3) with $a_{1}=2$, $a_{2}=0.02$, $a_{3}=0.5$, $a_{4}=2$,
$a_{5}=1.5$, $a_{6}=0,03$, $a_{7}=0.5$. The equilibrium point is:
$x_{0}=0.12,\quad y_{0}=0.99,\quad z_{0}=33,29$.

In the first case $k(s)=\delta(s-\tau)$, we obtain:
$\tau_{0}=0.8975032747$, $\omega_{0}=1.140149275$,
$g_{20}=-0.1032034629+0.8542910112i$,
$g_{11}=0.6157126005-0.5196432315i$,
$g_{02}=-0.1015556589+0.8738154402i$,
$g_{21}=-0.2181868792+0.9831278086i$,
$c_{1}(0)=-0.1197163151+0.2827563802i$, $\mu_{2}=0.2857035314$,
$\beta_{2}=-0.2394326302$, $T_2=0.1030243826$.

Because $\mu_{2}>0$, the Hopf bifurcation is supercritical for
$\tau>\tau_{0}$; as $\beta_{2}<0$ the bifurcating periodic
solution is asymptotically orbitally stable; as $T_{2}>0$ the
period increases. Here are the computer simulations. Fig.1
represents time vs uninfected cells $(t,x(t))$, Fig.2 time vs
infected cells $(t,y(t))$, Fig.3 time vs pathogens in blood
$(t,z(t))$, Fig.4 pathogens vs infected cells $(z(t), x(t))$,
Fig.5 pathogens vs infected cells $(z(t), y(t))$, Fig.6 uninfected
cells vs infected cells $(x(t),y(t))$.

\begin{center}\begin{tabular}{ccc}
\includegraphics[width=6cm]{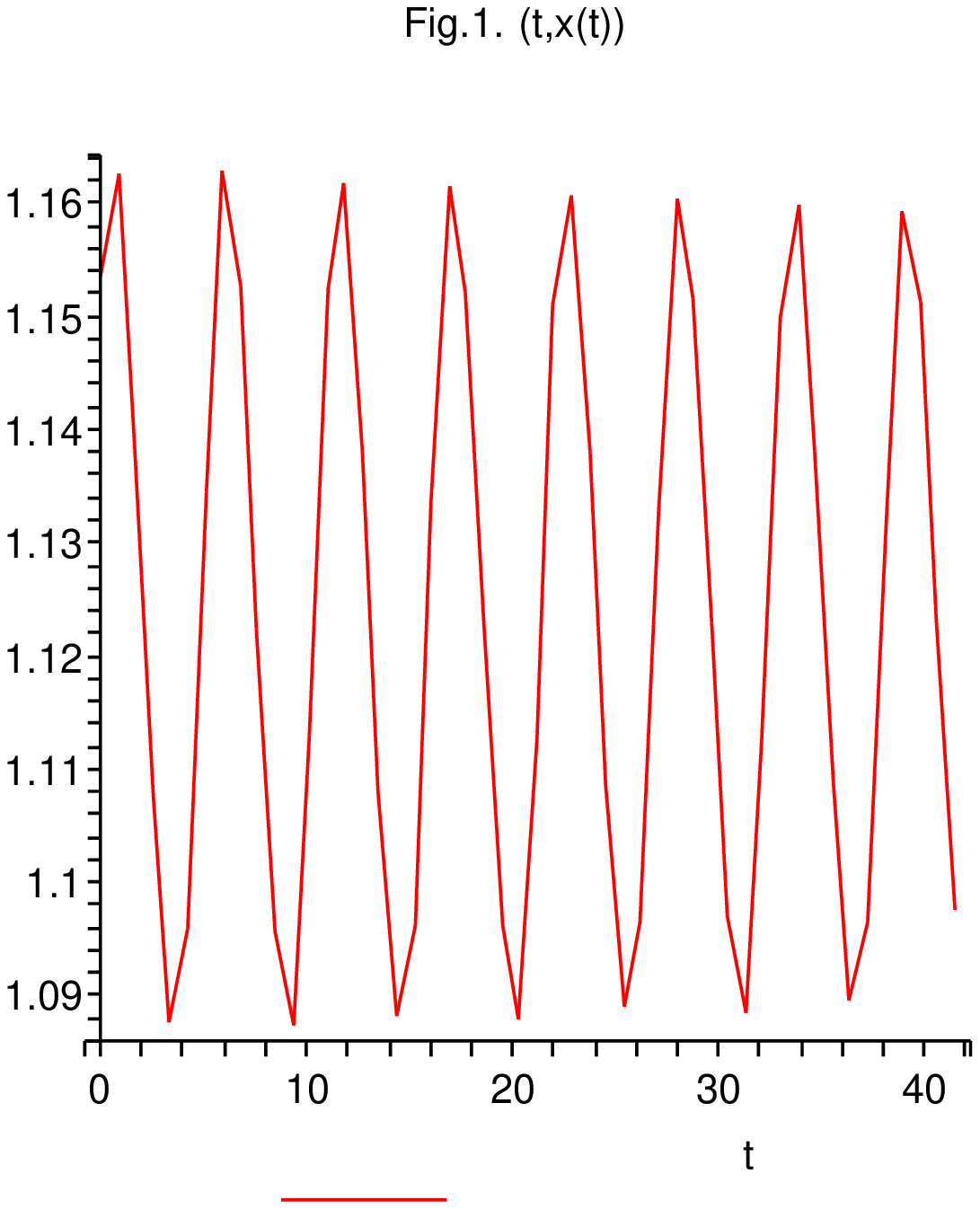} & &
\includegraphics[width=6cm]{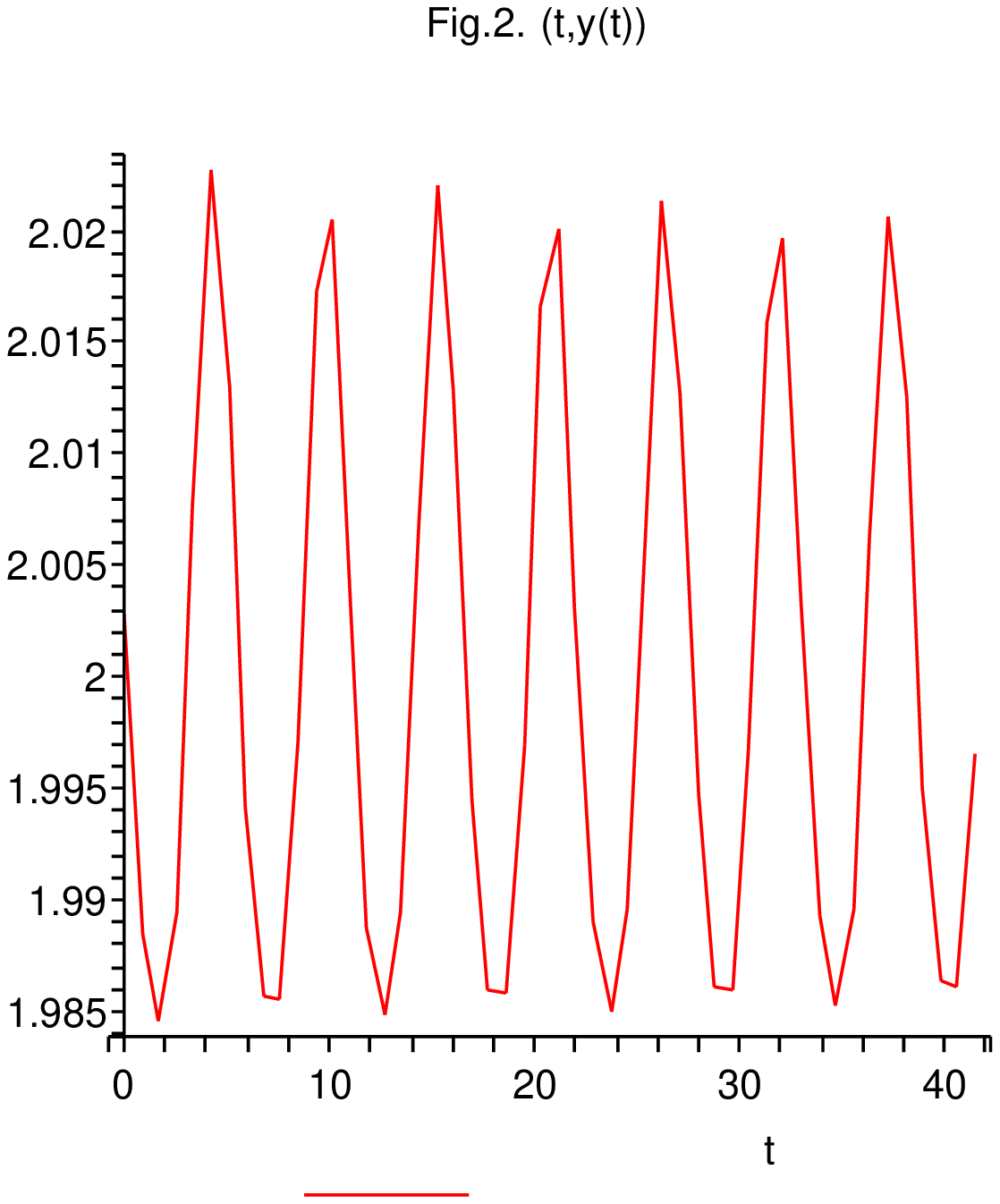}
\end{tabular}
\end{center}

\begin{center}
\begin{tabular}{ccc}
\includegraphics[width=6cm]{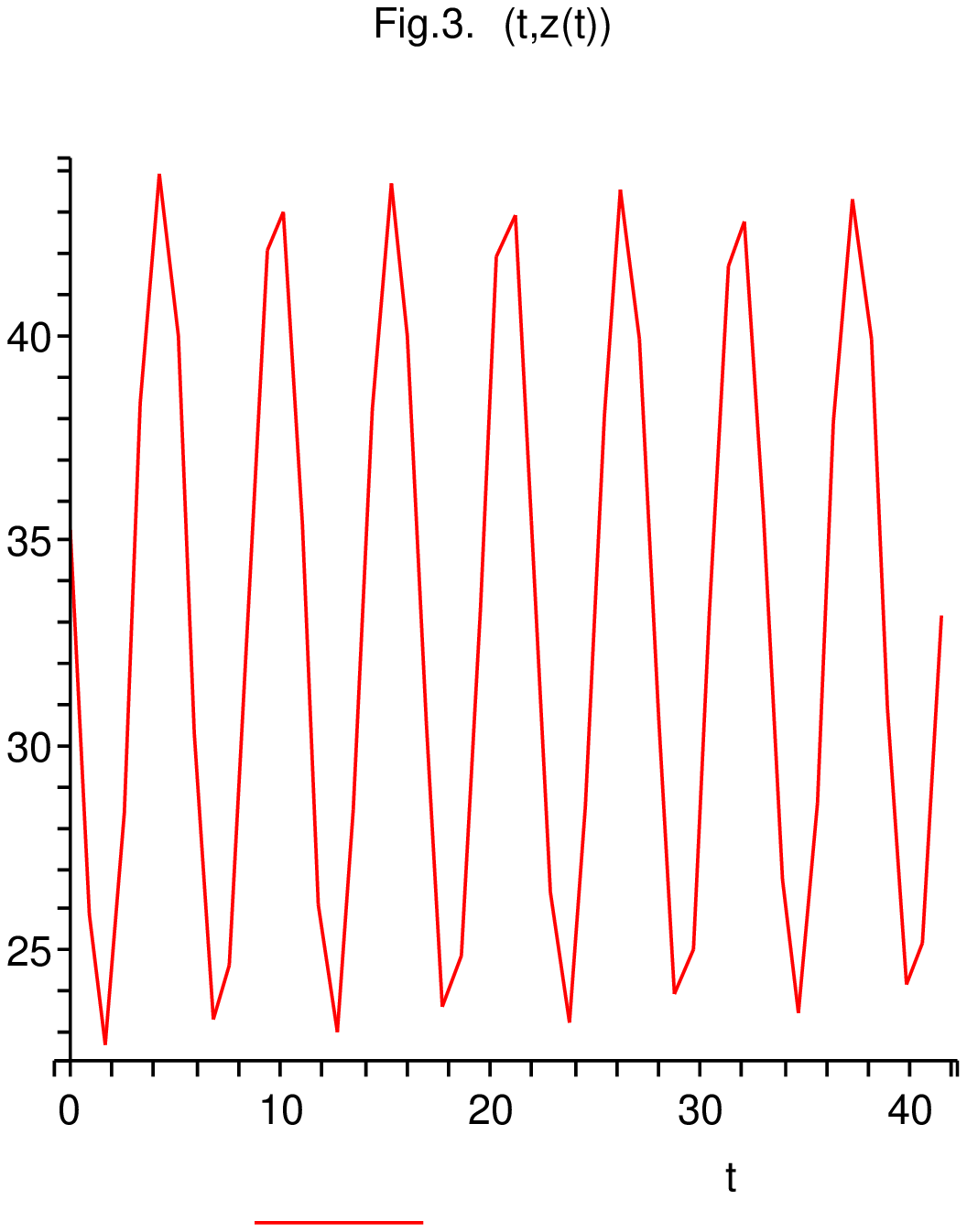} & &
\includegraphics[width=6cm]{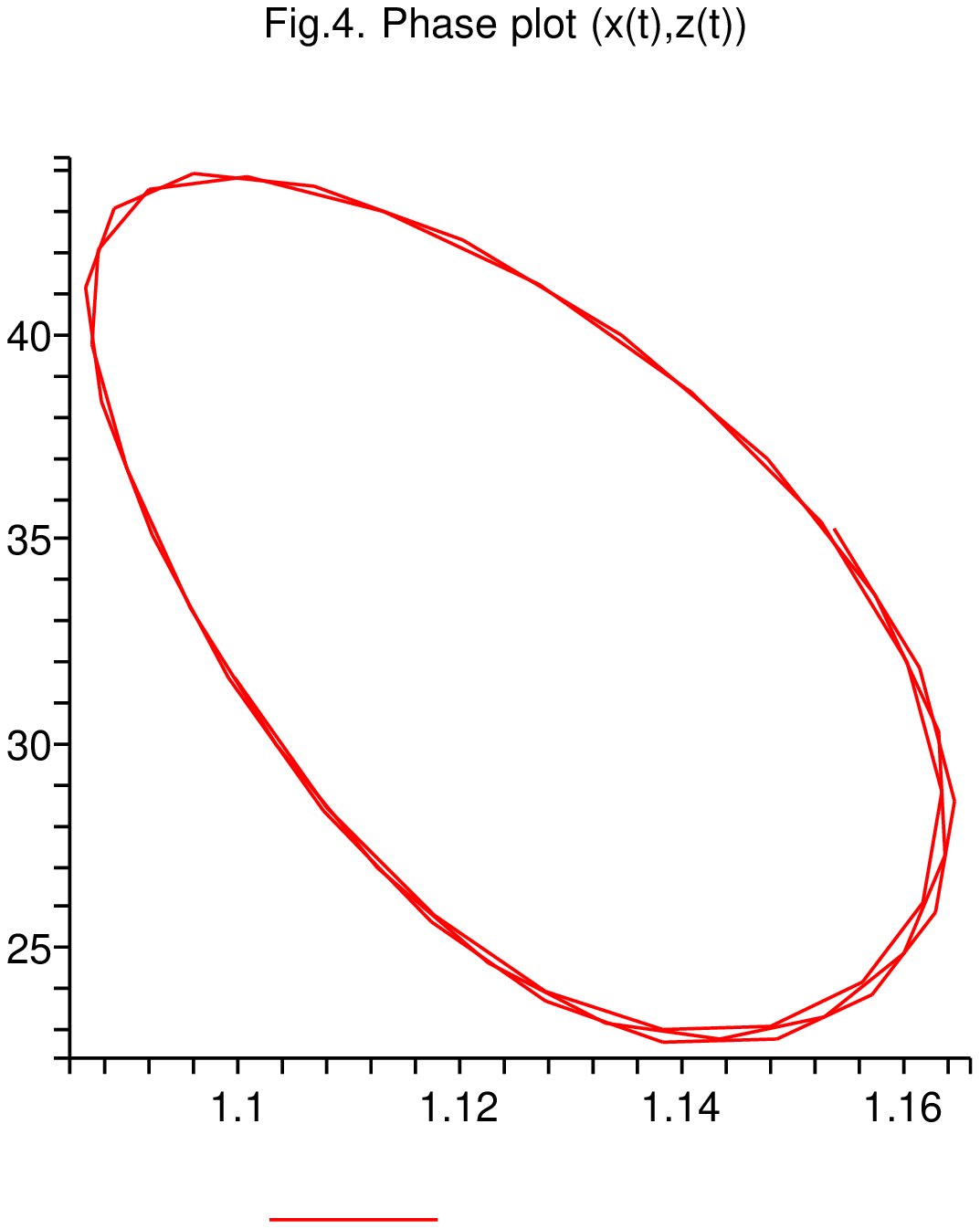}
\end{tabular}
\end{center}

\begin{center}
\begin{tabular}{ccc}
\includegraphics[width=6cm]{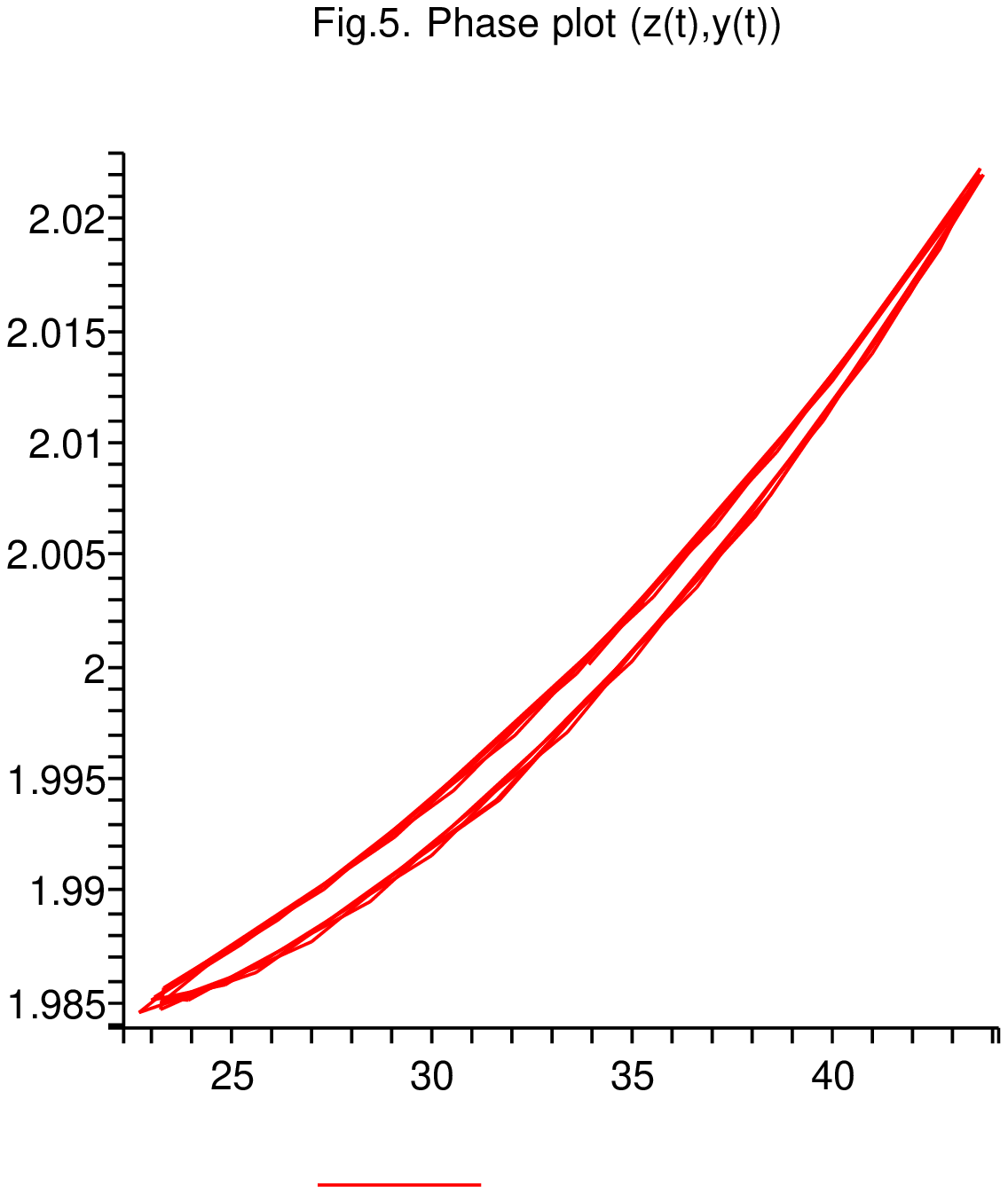}& &
\includegraphics[width=6cm]{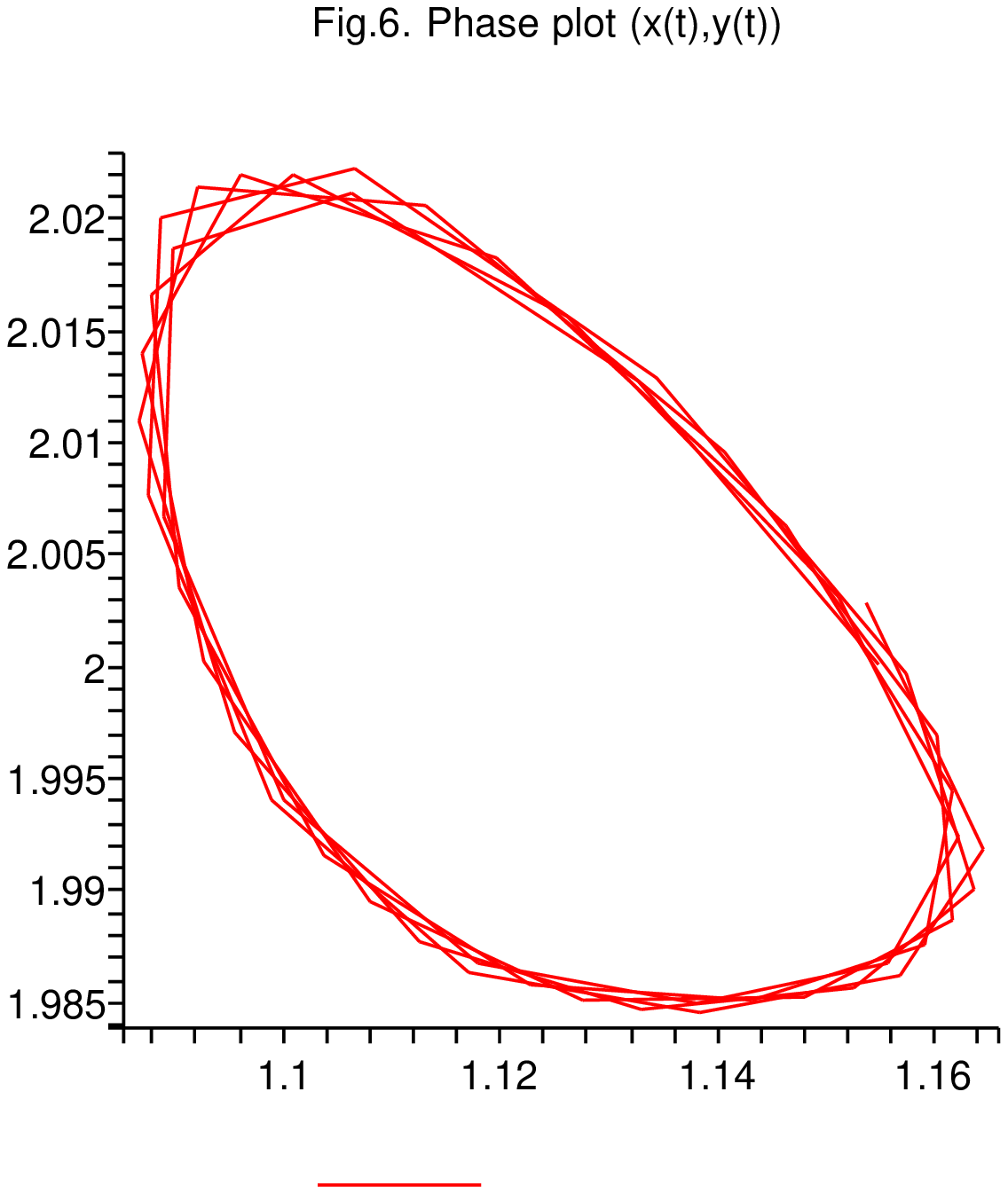}
\end{tabular}
\end{center}

In the second case $k(s)=qe^{-qs}$, we obtain: $q_0=0.1881852832$,
$\omega_{0}=0.1872904846$, $g_{20}=-0.1146541958-0.3262916164i$,
$g_{11}=0.2100798781-0.2554625514i$,
$g_{02}=-0.1034928343+0.5915640328i$,
$g_{21}=-0.1382064847-0.2196567760i$,
$c_{1}(0)=-0.6909931661-0.1103278787i$, $\mu_{2}=-0.4372331513$,
$\beta_{2}=-0.1381986332$, $T_2=0.8650529366$.

Because $\mu_{2}<0$, the Hopf bifurcation is subcritical for
$q>q_{0}$; as $\beta_{2}<0$ the bifurcating periodic solution is
asymptotically orbitally stable; as $T_{2}>0$ the period increase.
Here are the computer simulations: Fig.7 represents time vs
uninfected cells $(t,x(t))$, Fig.8 time vs infected cells
$(t,y(t))$, Fig.9 time vs pathogens in blood $(t,z(t))$, Fig.10
pathogens vs infected cells $(z(t), x(t))$, Fig.11 pathogens vs
infected cells $(z(t), y(t))$, Fig.12 uninfected cells vs infected
cells $(x(t),y(t))$.

\begin{center}\begin{tabular}{ccc}
\includegraphics[width=6cm]{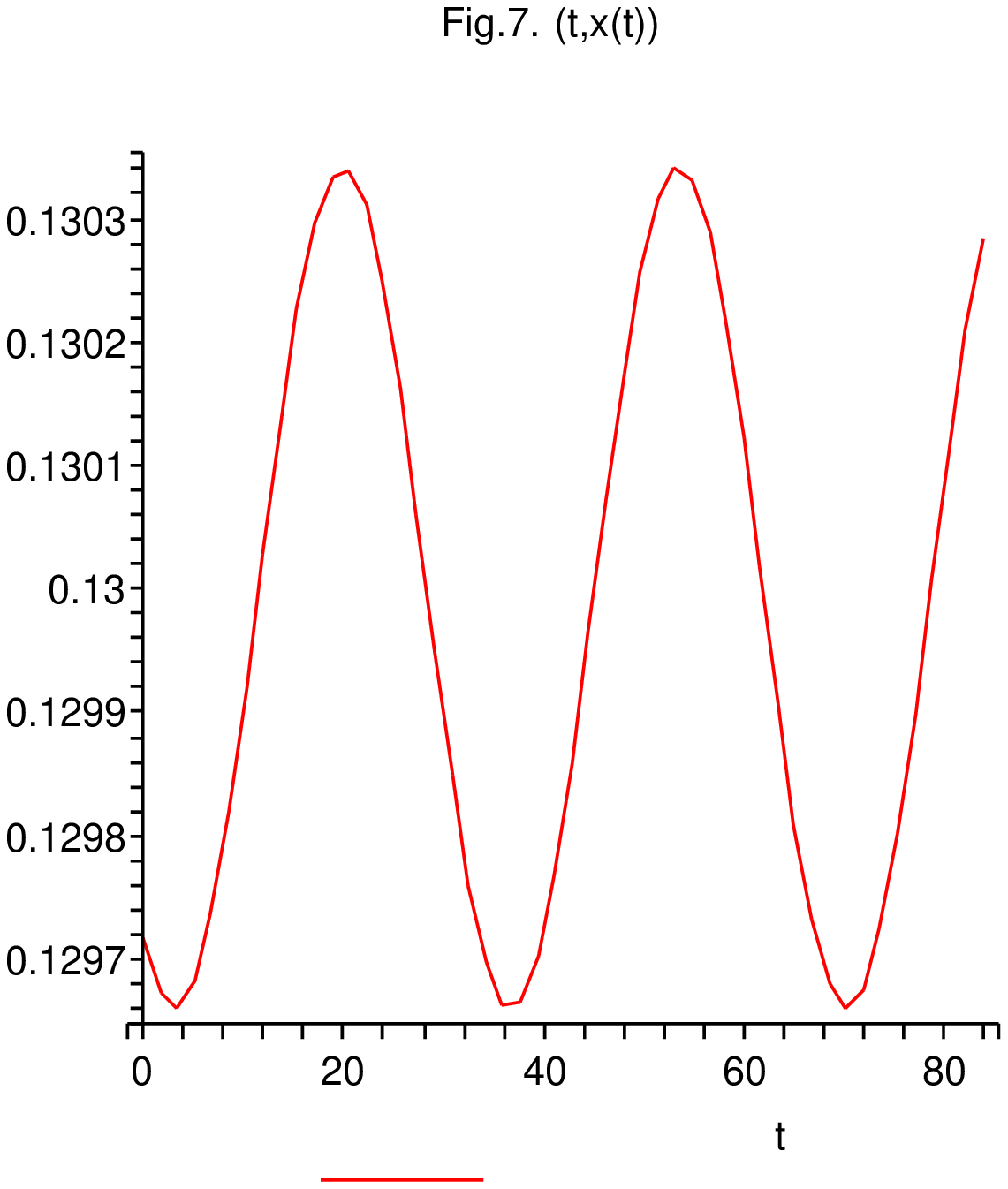} & &
\includegraphics[width=6cm]{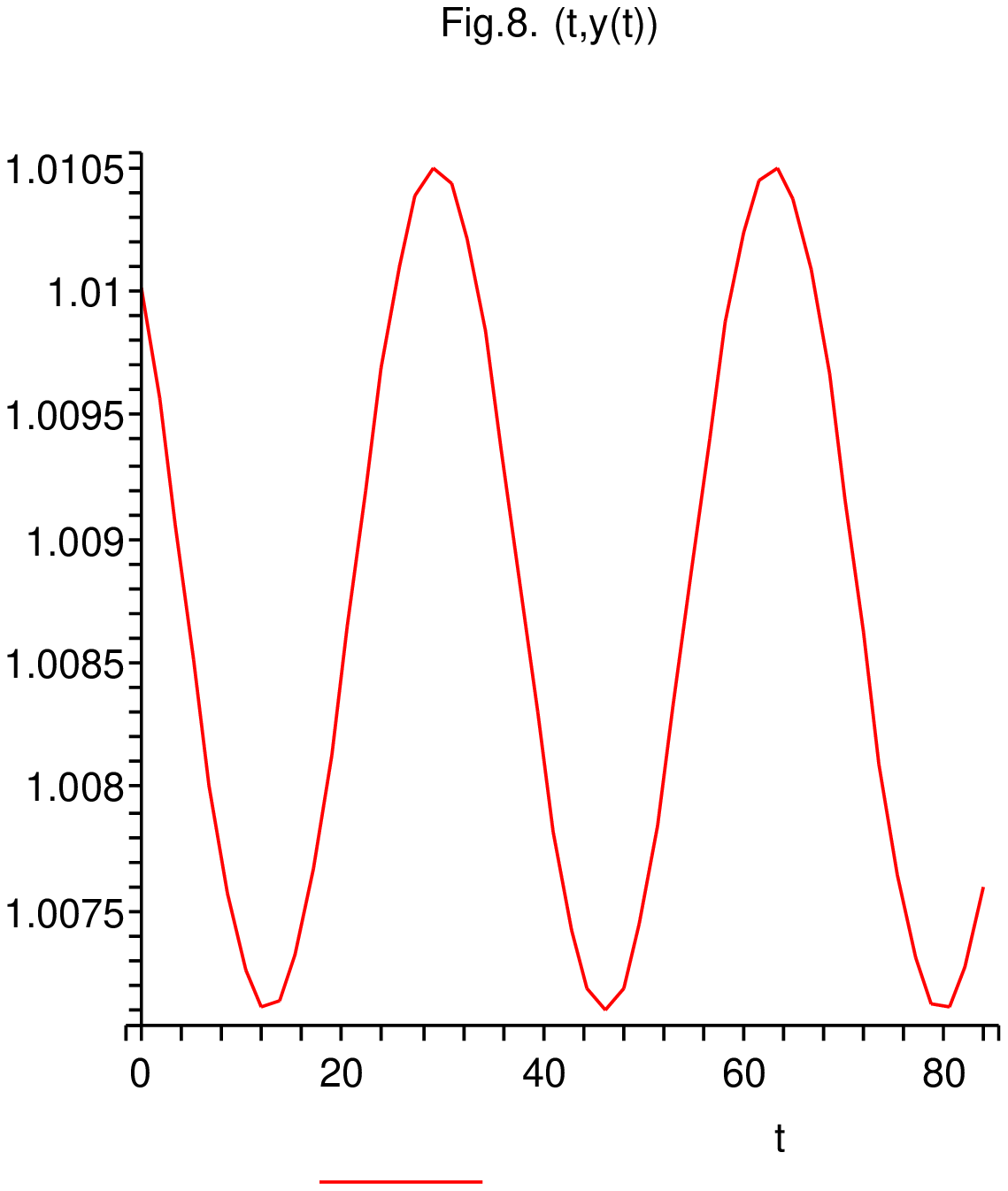}
\end{tabular}
\end{center}

\begin{center}
\begin{tabular}{ccc}
\includegraphics[width=6cm]{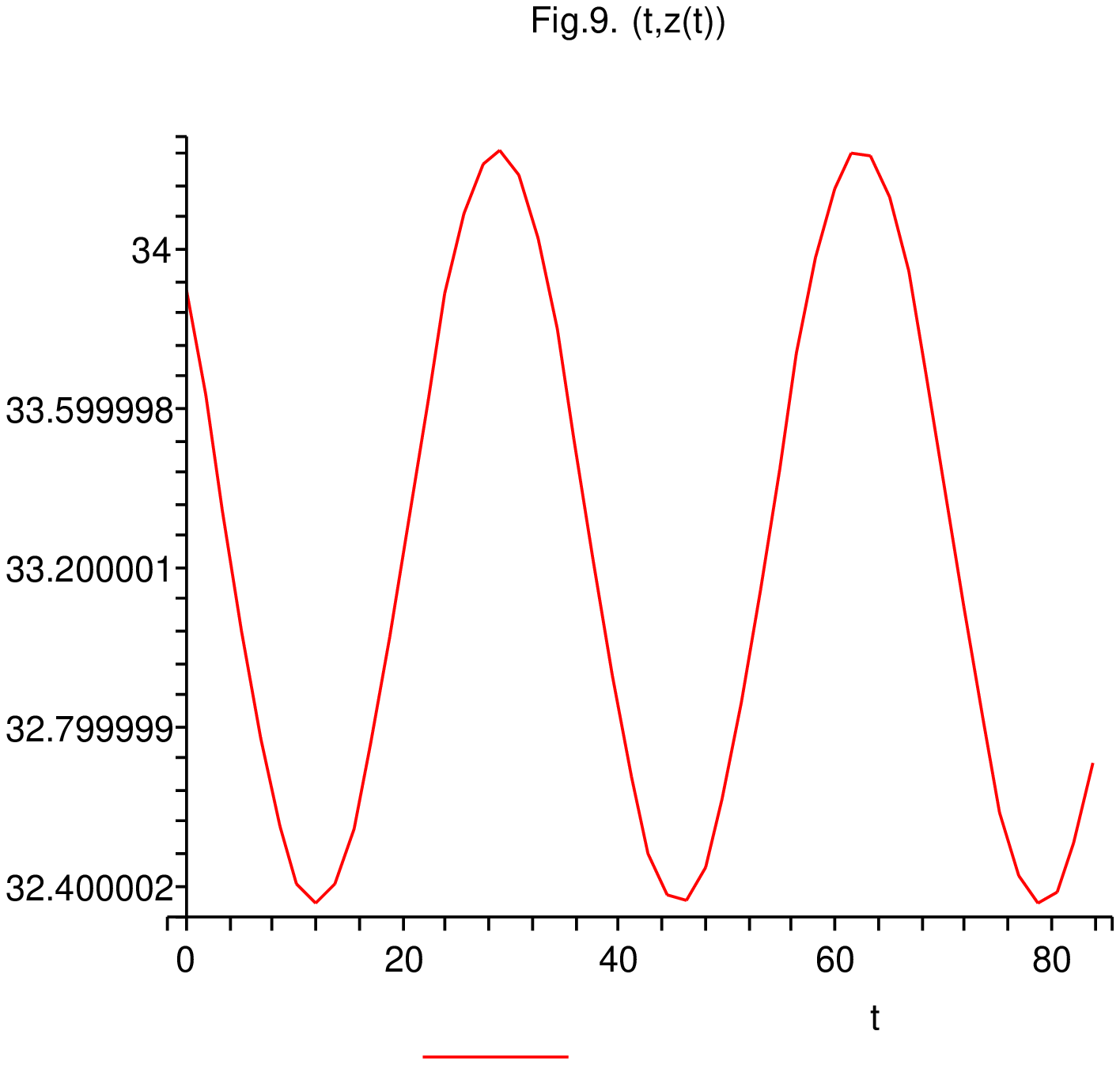} & &
\includegraphics[width=6cm]{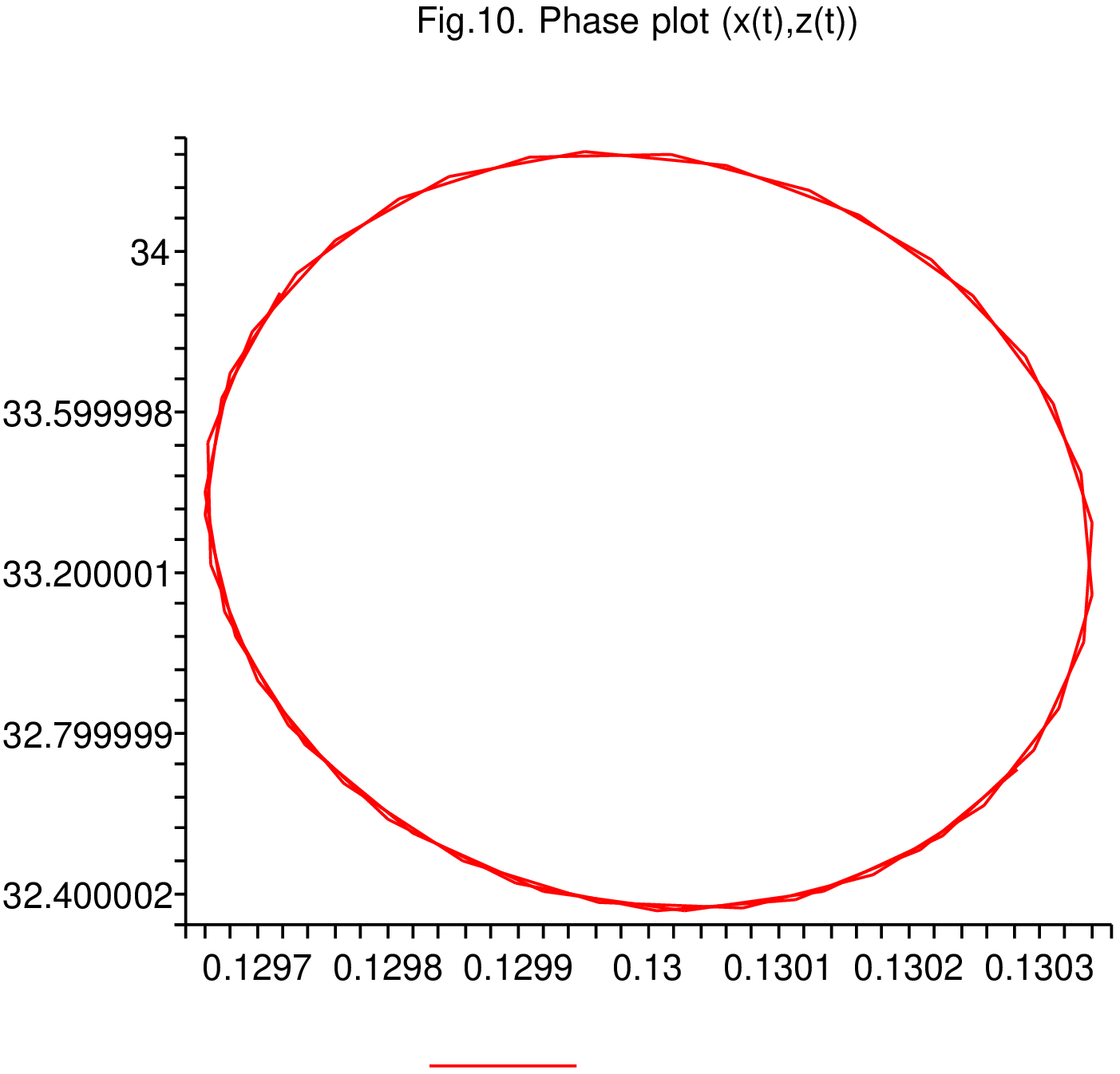}
\end{tabular}
\end{center}

\begin{center}
\begin{tabular}{ccc}
\includegraphics[width=6cm]{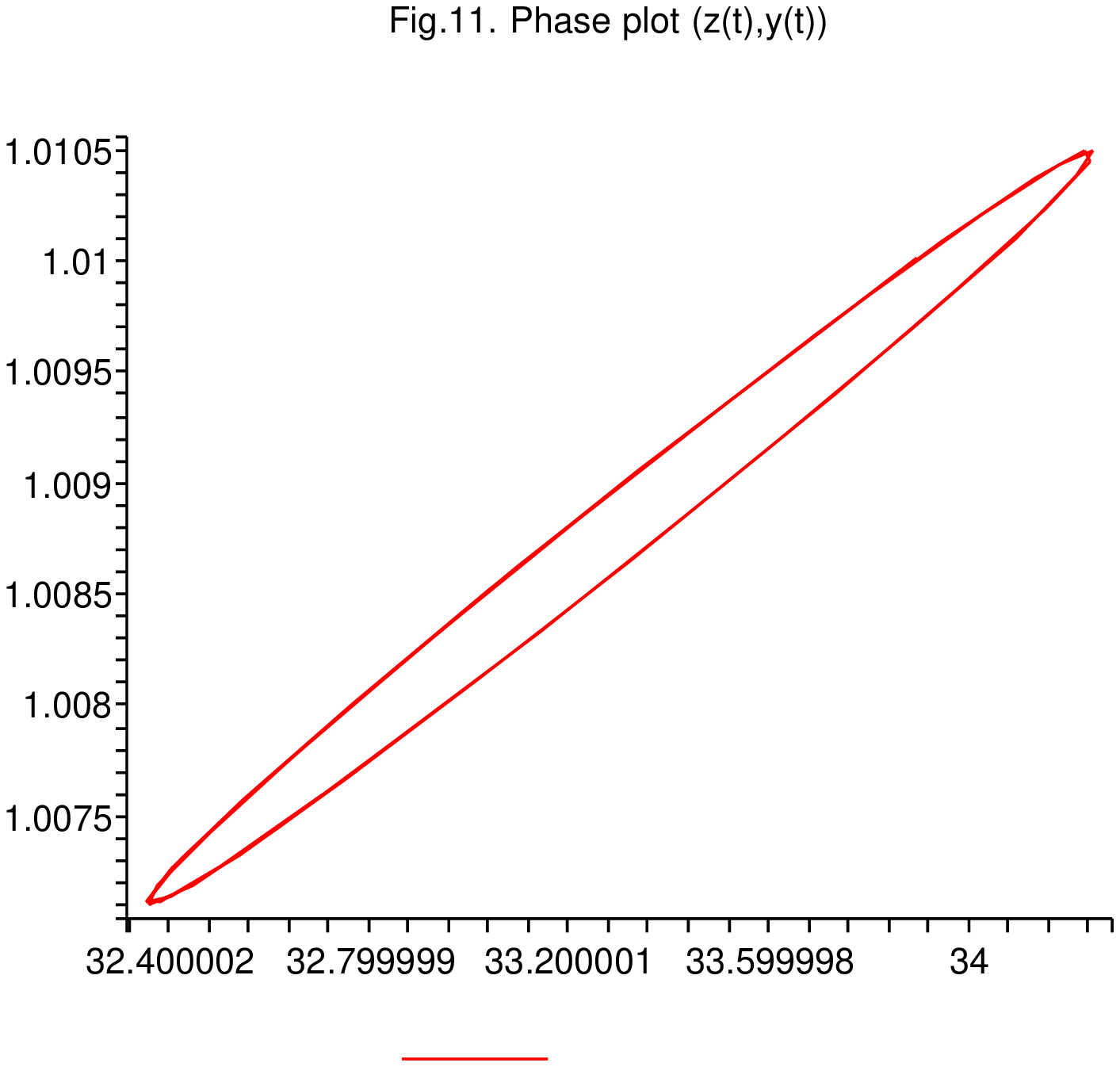}& &
\includegraphics[width=6cm]{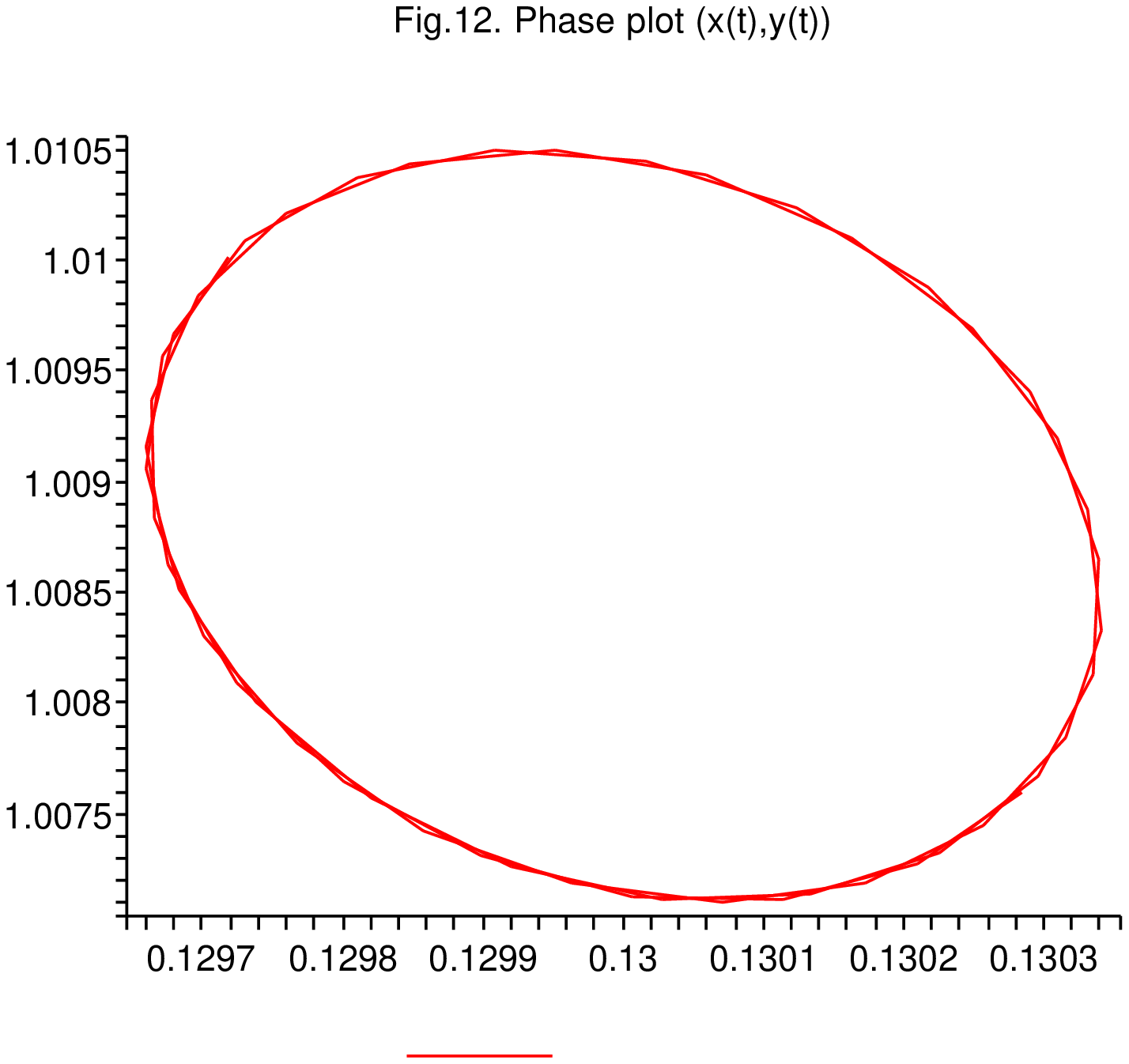}
\end{tabular}
\end{center}

\section{Conclusion}
In this paper we introduce a model which describes infectious
diseases and malaria infection with delay kernel. By using the
average time delay as  a parameter, it has been proved that the
Hopf bifurcation occurs when this parameter passes through a
critical value. A similarly  study will be made for the models
which describe infectious diseases, who takes into account immune
response against pathogens and the effect of involvement. This
study will be made in a future work.

\end{document}